\documentclass{amsart}
\usepackage[british]{babel}
\usepackage{mathrsfs}

\newtheorem{theorem}{Theorem}[section]
\newtheorem{lemma}[theorem]{Lemma}
\newtheorem{proposition}[theorem]{Proposition}

\theoremstyle{definition}
\newtheorem{definition}[theorem]{Definition}

\newtheorem{remark}[theorem]{Remark}

\newcommand{\sect}{\ensuremath{\S}~}

\newcommand{\field}[1]{\ensuremath{\mathbf{#1}}}
\newcommand{\Z}{\field{Z}}

\newcommand{\st}{\; {\big\vert} \;}
\newcommand{\Pres}[2]{\big\langle#1\st#2\big\rangle}
\newcommand{\subgroup}[1]{\langle#1\rangle}
\newcommand{\ab}[1]{\ensuremath{#1_{\text{ab}}}}
\newcommand{\chid}{\vec{\chi}_2}

\newcommand{\lift}[1]{\ensuremath{\widetilde{#1}}}

\newcommand{\ZG}{\ensuremath{\Z\Gamma}}
\newcommand{\FPnn}[1]{\ensuremath{\mathrm{FP}_{#1}}}

\DeclareMathOperator{\rk}{rk}
\DeclareMathOperator{\im}{im}
\DeclareMathOperator{\defic}{def}
\DeclareMathOperator{\adef}{adef}

\begin{document}
\title[Deficiency and abelianized deficiency]{Deficiency and abelianized
deficiency of some virtually free groups}

\author{Martin R.~Bridson}
\address{Dept.\ of Mathematics, Imperial College, London SW7 2AZ.}
\email{m.bridson@imperial.ac.uk}
\thanks{}

\author{Michael Tweedale}
\address{Dept.\ of Mathematics, University of Bristol, Bristol BS8 1TW.}
\email{m.tweedale@bristol.ac.uk}
\thanks{This work was supported in part by grants from the EPSRC. The first
author is also supported  by a Royal Society Wolfson Research Merit Award.}

\subjclass[2000]{20F05 (primary), 57M20 (secondary)}

\date{2nd November 2006}
\begin{abstract}
  Let $Q_m$ be the HNN~extension of $\Z/m\times\Z/m$ where the stable letter
  conjugates the first factor to the second. We explore small presentations of
  the groups $\Gamma_{m,n}=Q_m\ast Q_n$. We show that for certain choices of
  $(m,n)$, for example $(2,3)$, the group $\Gamma_{m,n}$ has a relation gap
  unless it admits a presentation with at most~$3$ defining relations, and we
  establish restrictions on the possible form of such a presentation. We then
  associate to each $(m,n)$ a $3$-complex with $16$~cells. This $3$-complex is
  a counterexample to the $D(2)$~conjecture if $\Gamma_{m,n}$ has a relation
  gap.
\end{abstract}
\maketitle\let\languagename\relax\sloppy

\section{Introduction}\label{sec:intro}
\noindent Given a finite presentation $F/R$ for a group $\Gamma$, the action
of the free group $F$ by conjugation on $R$ induces an action of $\Gamma$ on
the abelian group $\ab{R}=R/[R,R]$. It is obvious that the rank of $\ab{R}$ as
a $\ZG$-module is at most the number of elements needed to generate $R$ as a
normal subgroup of $F$; it therefore serves as a lower bound on the minimal
number of relators needed to present $\Gamma$ on the given generators. This
lower bound seems extremely crude: one can hardly believe that it will be
sharp in general. And yet, despite sustained attack over many years, not a
single example has emerged to lend substance to this intuition.

A finite presentation $F/R$ whose \emph{relation module} $\ab{R}$ has rank
strictly smaller than the number of elements required to generated $R$ as a
normal subgroup of $F$ is said to have a \emph{relation gap}; the
\emph{relation gap problem} is to determine whether or not there exists a
presentation with a relation gap. It belongs to a circle of famous and
notoriously hard problems concerning the homotopy properties of finite
$2$-complexes.  For example, it is closely related to the
\emph{$D(2)$~conjecture}.

In~\cite{Wall}, C.~T.~C.~Wall established that homological invariants are
sufficient to determine whether a CW~complex has the homotopy type of an
$n$-dimensional CW~complex, except possibly when $n=2$. Recall that a space
$X$ has \emph{property~$D(n)$} if $H_i(\lift{X})=0$ for $i>n$, where
$\lift{X}$ is the universal cover of $X$, and if in addition
$H^{n+1}(X;\mathscr{M})=0$ for all local coefficient systems $\mathscr{M}$ on
$X$. Wall's result is that if $n\ne2$, then a (finite) CW~complex has the
homotopy type of a (finite) $n$-complex if and only if it has property~$D(n)$.
The assertion that a finite $3$-complex has the homotopy type of a finite
$2$-complex if and only if it has property~$D(2)$ has become known as
\emph{the $D(2)$~conjecture}. A theorem of M.~Dyer (unpublished; a proof can
be found in~\cite{HarlSome}) states that if a group $\Gamma$ with
$H^3(\Gamma;\ZG)=0$ has a presentation with a relation gap that \emph{realizes
the deficiency} of the group (see \sect\ref{sec:background} for a discussion
of deficiency), then the $D(2)$~conjecture is false.

The purpose of this note is to propose a new collection of candidates for
presentations with relation gaps. The prime merit of these examples is that
one can give a short, transparent and natural proof that the relation modules
of the obvious presentations can be generated by one fewer element than one
would expect. (Using similar ideas one can create many other examples, but
we have resisted the temptation to present these because we do not want to
obscure the main idea: the straightforwardness of our examples is what we find
most attractive about them.)

Let $Q_m$ be the HNN~extension of $\Z/m\times\Z/m$ where the stable letter
conjugates the first factor to the second. The obvious presentation of
$Q_m\ast Q_n$ has four generators and four relations, but we shall see that if
$m$ and $n$ satisfy a coprimeness condition then the relation module
associated to this presentation requires only three generators. The argument
extends easily to arbitrarily many free factors $Q_{m_1}\ast\cdots\ast
Q_{m_r}$, yielding presentations where the expected relation gap is $r-1$.
Using a result of J.~Howie~\cite{Howie} on one-relator products of locally
indicable groups, we obtain restrictions on putative 3-relator presentations of
$Q_m\ast Q_n$ (Proposition~\ref{prop:nopres}), but as yet we have been unable
to prove that such presentations do not exist. The difficulty of doing so is
discussed in \sect\ref{sec:lower}.

In \sect\ref{sec:d2} we pursue the line of attack on the $D(2)$~conjecture via
relation gaps discussed above. We follow Harlander's construction to give an
explicit description of a $3$-complex with only $16$~cells that looks
homologically like a $2$-complex, in the sense that it possesses Wall's
property~$D(2)$, but that does not have the homotopy type of a finite
$2$-complex if $Q_m\ast Q_n$ has a relation gap.

Groups similar to $Q_m\ast Q_n$ have been studied previously for their
interesting presentation theory. For example, the groups used
in~\cite{HogLustigMetzler} to show that deficiency is not additive under the
operation of free product were of the form
$(\Z/m\times\Z/m')*(\Z/n\times\Z/n')$. The groups
$G_{m,n}=(\Z/m\times\Z)*(\Z/n\times\Z)$ were studied in connection with
efficiency by D.~Epstein~\cite{EpsteinFinite}, and have come to the fore again
recently in the work of K.~Gruenberg and P.~Linnell~\cite{GruenbergLinnell}.
This last work, which is more sophisticated than ours, focuses on the
presentation theory of free products of finite groups, but also contains a
proof that the relation module associated to the obvious $4$-generator,
$4$-relation presentation of $G_{m,n}$ requires only three generators.

We should mention that the obvious extension of the relation gap problem to
finitely generated, rather than just finitely presented, groups is also an
interesting problem, and has been resolved by M.~Bestvina and N.~Brady:
in~\cite{BestvinaBrady}, they construct finitely generated groups that are not
finitely presented but have finitely generated relation modules. These groups
therefore have `infinite relation gaps'. Naturally, one thinks of trying to
build on these examples to find finitely presented groups with relation gaps:
we take up this approach in~\cite{BTArtin} (see also~\cite{HowieBB}).

\section{Deficiency and abelianized deficiency}\label{sec:background}
\noindent In this section, we assemble some basic definitions and well-known
results for later reference. We write $d(\Gamma)$ for the minimum number of
elements needed to generate a group $\Gamma$. If $Q$ is a group acting on
$\Gamma$ then we write $d_Q(\Gamma)$ for the minimum number of $Q$-orbits
needed to generate $\Gamma$.

Let $\Gamma$ be a finitely presented group. The \emph{deficiency} of a finite
presentation $F/R$ of $\Gamma$ is $d_F(R)-d(F)$, where $F$ operates on its
normal subgroup $R$ by conjugation. (Some authors' definition of deficiency
differs from ours by a sign.)

The action of $F$ on $R$ induces by passage to the quotient an action of
$\Gamma$ on the abelianization $\ab{R}$ of $R$, which makes $\ab{R}$ into a
$\ZG$-module, called the \emph{relation module} of the presentation. The
\emph{abelianized deficiency} of the presentation is
$d_{\Gamma}(\ab{R})-d(F)$; this invariant was first studied by K.~Gruenberg,
under the name of \emph{abelianized defect}.

\begin{lemma}\label{lem:relmodexact}
  If $F$ is free of rank~$d$, then there is an exact sequence of $\ZG$-modules
  \[
  0\to\ab{R}\to(\ZG)^d\to\ZG\to\Z\to0.
  \]
\end{lemma}
\begin{proof}
  We identify $F$ with the fundamental group of a graph $X$ that has one
  vertex and $d$ edges. The regular covering $\hat{X}$ of $X$ corresponding to
  the subgroup $R\subset F$ is the Cayley graph of $\Gamma$. The exact
  sequence in the statement of the lemma is obtained from the cellular chain
  complex of this covering by inserting the first homology group
  $H_1(\hat{X})=\ab{R}$ on the left as the kernel of the first boundary map.
\end{proof}

\begin{lemma}\label{lem:defdef}
  The deficiency of any finite presentation of $\Gamma$ is bounded below by
  the abelianized deficiency, and this in turn is bounded below by
  $d(H_2(\Gamma))-\rk(H_1(\Gamma))$, where $\rk$ is torsion-free rank.
\end{lemma}
\begin{proof}
  The first part is clear. For the second, let $C_*$ be the cellular complex
  of the universal cover of a presentation $2$-complex of $\Gamma$: this is
  a partial resolution of $\Z$ by free $\ZG$-modules. Apply the functor
  $-\otimes_{\ZG}\Z$ to $C_*$. If $C_i$ has rank $r_i$, then the resulting
  complex has a free abelian group of the same rank $r_i$ in degree $i$;
  moreover, the homology groups of the new complex are $H_i(\Gamma)$ in
  degrees $i=0,1$, and $\ker\partial_2$ in degree $2$. Since $H_2(\Gamma)$
  is a quotient of this kernel, one has
  $d(H_2(\Gamma))\le\rk(\ker\partial_2)\le r_2$, and the result follows.
\end{proof}

\noindent Therefore we can define:

\begin{definition}
  The \emph{deficiency} $\defic(\Gamma)$ (resp.\ \emph{abelianized deficiency}
  $\adef(\Gamma)$) of $\Gamma$ is the infimum of the deficiencies (resp.\
  abelianized deficiencies) of the finite presentations of $\Gamma$.
\end{definition}

\begin{remark}
  Obviously, if $\Gamma$ has a presentation of deficiency
  $d(H_2(\Gamma))-\rk(H_1(\Gamma))$, then by Lemma~\ref{lem:defdef} this
  presentation realizes the deficiency of the group. In this case, $\Gamma$ is
  said to be \emph{efficient}. One knows that inefficient groups exist: R.~Swan
  constructed finite examples in~\cite{Swan}, and much later
  M.~Lustig~\cite{LustigEff} produced the first torsion-free examples. Further
  examples are given by superperfect groups that are not fundamental groups of
  homology $4$-spheres (see~\cite{HausmannWeinberger} and \cite{Hillman}).
\end{remark}

\section{The examples}\label{sec:examples}
\noindent Given letters $x$ and $t$, let $\rho_n=\rho_n(x,t)$ be the word
\[
\rho_n(x,t)=(txt^{-1})x(txt^{-1})^{-1}x^{-n-1}.
\]
We will consider the groups $Q_n = \Pres{x,t}{\rho_n,x^{n}}$.

\begin{lemma}\label{lem:hnnpres}
  $Q_n$ is isomorphic to the HNN~extension
  \[
  (\Z/n\times\Z/n)\ast_{\phi},
  \]
  where $\phi$ maps the first factor isomorphically to the second.
\end{lemma}
\begin{proof}
  We apply Tietze moves to the given presentation of $Q_n$, first adding a
  superfluous generator, and then adding a redundant relator:
  \begin{align*}
    Q_n &= \Pres{x,t}{[(txt^{-1}),x],x^{n}}\\
    &= \Pres{x,b,t}{b=txt^{-1},[b,x],x^{n}}\\
    &= \Pres{x,b,t}{b=txt^{-1},[b,x],x^{n},b^{n}}.
  \end{align*}
\end{proof}

\subsection{Generators for the relation module}\label{sec:relmodgens}
Let $q_n=(n+1)^{n}-1$ and $c_n=n\,q_n$.

\begin{lemma}\label{lem:cn}
  In $\Pres{x,t}{\rho_n}$, one has the equality
  $[(txt^{-1})^{n},x^{n}]=x^{c_n}$.
\end{lemma}
\begin{proof}
  Conjugating $x$ by $txt^{-1}$ in our group $n$ times, we have
  \[
  (txt^{-1})^{n}x(txt^{-1})^{-n}=x^{q_n+1}.
  \]
  Raising this to the power $n$ gives
  \[
  (txt^{-1})^{n}x^{n}(txt^{-1})^{-n}=x^{n(q_n+1)},
  \]
  i.e.
  \[
  [(txt^{-1})^{n},x^{n}]=x^{c_n}.
  \]
\end{proof}

\begin{proposition}\label{prop:relmodgens}
  Suppose that $(q_m,q_n)=1$, and let
  \begin{align*}
    \Gamma_{m,n}&=Q_m*Q_n\\
    &=\Pres{x_m,t_m,x_n,t_n}{\rho_m(x_m,t_m),\rho_n(x_n,t_n),x_m^{m},x_n^{n}}.
  \end{align*}
  Let $\ab{R}$ be the relation module of this presentation. Then $\ab{R}$ is
  generated as a $\Z\Gamma_{m,n}$-module by the images of $\rho_m$, $\rho_n$
  and $x_m^{m}x_n^{n}$.
\end{proposition}
\begin{proof}
  Let $S$ denote the quotient of $\ab{R}$ by the $\ZG$-submodule generated
  by $\rho_m$ and $\rho_n$, and let $\pi:R\to S$ be the obvious surjection.
  It is clear that $S$ is generated as a $\ZG$-module by the images of
  $x_m^{m}$ and $x_n^{n}$; it will be sufficient for us to show that in fact
  $S$ is a cyclic $\ZG$-module, generated by the image of $x_m^{m}x_n^{n}$.

  Since $x_m^{m}=(t_mx_mt_m^{-1})^{m}=1$ in $Q_m$, both $x_m^{m}$ and
  $(t_mx_mt_m^{-1})^{m}$ lie in $R$, and hence their commutator lies in
  $[R,R]$. On the other hand, by the previous lemma $x_m^{c_m}$ is equal to
  this commutator modulo $\rho_m$, so $\pi(x_m^{c_m})=1$. Thus the order of
  $\pi(x_m^{m})$ in $S$ divides $q_m=c_m/m$. Similarly, the order of
  $\pi(x_n^{n})$ in $S$ divides $q_n$. Since $q_m$ and $q_n$ are coprime, it
  follows that $S$, which is generated by $\pi(x_m^{m})$ and $\pi(x_n^{n})$,
  is actually generated by $\pi(x_m^{m})\pi(x_n^{n})$ alone, as claimed.
\end{proof}

An entirely similar argument yields:

\begin{proposition}\label{prop:morefactors}
  If $(q_{m_i},q_{m_j})=1$ for $1\le i< j\le r$, then the relation module of
  \begin{align*}
    \Gamma&=Q_{m_1}\ast\cdots\ast Q_{m_r}\\
    &=\Pres{x_{m_1},t_{m_1},\ldots,x_{m_r},t_{m_r}}{\rho_{m_1}(x_{m_1},t_{m_1}),\,
    x_{m_1}^{m_1},\ldots,\rho_{m_r}(x_{m_r},t_{m_r}),\, x_{m_r}^{m_r}}
  \end{align*} is
  generated as a $\Z\Gamma$-module by the images of
  $\rho_{m_1},\ldots,\rho_{m_r}$ and $x_{m_1}^{m_1}\ldots x_{m_r}^{m_r}$.
\end{proposition}

\subsection{An observation on putative 3-relator presentations}\label{sec:howie}
Although it certainly does not approach a proof that our groups have relation
gaps, the result of this subsection restricts the nature of possible
presentations.

Recall that a group is \emph{locally indicable} if each of its non-trivial
finitely generated subgroups has the infinite cyclic group as a homomorphic
image. We will apply a result about these groups due to J.~Howie.

\begin{theorem}[{\cite[Theorem~4.2]{Howie}}]\label{th:howie}
  Let $A$ and $B$ be locally indicable groups, and let $G$ be the quotient of
  $A*B$ by the normal closure of a single element $r$, not conjugate to an
  element of $A$ or $B$. The following are equivalent:
  \begin{enumerate}
    \item\label{item:howiea} $G$ is locally indicable;
    \item\label{item:howieb} $G$ is torsion-free;
    \item\label{item:howiec} $r$ is not a proper power in $A*B$.
  \end{enumerate}
\end{theorem}

\noindent We call such a group $G$ a \emph{one-relator product} of $A$ and
$B$.

\begin{proposition}\label{prop:nopres}
  In the notation of Proposition~\ref{prop:relmodgens}, the group
  $\Gamma_{m,n}$ does not admit a presentation of the form
  \[
  \Pres{x_m,t_m,x_n,t_n}{\rho_m,\rho_n,r},
  \]
  for any word $r\in\{x_m,t_m,x_n,t_n\}^*$.
\end{proposition}
\begin{proof}
  Consider $L=A*B$, where $A=\Pres{x_m,t_m}{\rho_m}$ and
  $B=\Pres{x_n,t_n}{\rho_n}$. This is a free product of torsion-free
  one-relator groups (in particular, of locally indicable groups). Suppose
  for a contradiction that $\Gamma_{m,n}$ is a one-relator quotient of $L$
  by $r$. Since $\Gamma_{m,n}$ has torsion, Theorem~\ref{th:howie} implies
  that either $r$ is conjugate in $L$ to an element of $A$ or $B$, or else
  $r=r_0^p$ for some $r_0\in L$, not itself a proper power. In the second
  case, the image of $r_0$ in $\Gamma_{m,n}$ (which we also denote by $r_0$)
  is a torsion element, and so is conjugate in $\Gamma_{m,n}$ to an element
  of $H_m=\subgroup{x_m,t_mx_mt_m^{-1}}\cong\Z/m\times\Z/m$ or
  $H_n=\subgroup{x_n,t_nx_nt_n^{-1}}\cong\Z/n\times\Z/n$: for the sake of
  argument, let us say of $H_m$.  Consider $M$, the one-relator quotient of
  $L$ by $r_0$. The quotient map $L\to M$ factors through $\Gamma_{m,n}$, so
  $M$ is isomorphic to $\Gamma_{m,n}$ quotiented by the single additional
  relator $r_0$. Since $r_0$ is conjugate in $\Gamma_{m,n}$ to an element
  $r_0'$ of $H_m$, it follows that $M$ is also isomorphic to $\Gamma_{m,n}$
  quotiented by $r_0'$.  But it is easy to see that in this last group, the
  image of $x_n$ is still a torsion element. Thus
  condition~\ref{item:howiec} of Theorem~\ref{th:howie} holds for $M$, but
  condition~\ref{item:howieb} fails, so we conclude that the hypothesis of
  the theorem cannot be satisfied: in other words, $r_0$ must be conjugate
  in $L$ to an element of $A$ or $B$.

  In either case, after conjugating $r$ if necessary, we can assume that $r$
  is contained in one of the factors, $A$ or $B$: let's say in $A$.  But
  this means that $\Gamma_{m,n}$ splits as a free product
  \[
  \Gamma_{m,n}=\Pres{x_m,t_m}{\rho_m,r}*\Pres{x_n,t_n}{\rho_n},
  \]
  which is absurd since $x_n$ has order~$n$ in $\Gamma_{m,n}$.
\end{proof}

\section{A proposed counterexample to the $D(2)$~conjecture}\label{sec:d2}
\noindent In this section we associate to each of our groups $\Gamma_{m,n}$ a
finite $3$-complex with only $16$~cells. If $\Gamma_{m,n}$ has a relation gap
then the corresponding 3-complex will be a counterexample to the
$D(2)$~conjecture.

\begin{theorem}[$\mathrm{Dyer}$]\label{th:dyer}
  Let $\Gamma$ be a group with $H^3(\Gamma,\Z\Gamma)=0$. If there is a
  presentation of $\Gamma$ that realizes the deficiency of the group and that
  has a relation gap, then the $D(2)$~conjecture is false.
\end{theorem}

A proof of this theorem, including a construction of a counterexample to the
$D(2)$~conjecture under the given hypotheses, has been given by Jens Harlander
in his fine survey article~\cite{HarlSome}. In this section we follow
Harlander's construction for our groups $\Gamma_{m,n}$.

Let $K$ be the presentation $2$-complex associated to our original
presentation of $\Gamma=\Gamma_{m,n}$,
\[
\Gamma=\Pres{x_m,t_m,x_n,t_n}{\rho_m(x_m,t_m),\rho_n(x_n,t_n),x_m^{m},x_n^{n}}.
\]
By Proposition~\ref{prop:relmodgens}, the relation module $\ab{R}$ associated
to this presentation can be generated by three elements, and so there is an
exact sequence of $\ZG$-modules
\begin{equation*}
  \tag{$*$}0\to N\to(\ZG)^3\to\ab{R}\to0.
\end{equation*}
On the other hand, the cellular chain complex of the universal cover of $K$ is
a complex of $\ZG$-modules $0\to
(\ZG)^4\xrightarrow{\partial_2}(\ZG)^4\xrightarrow{\partial_1}\ZG\to\Z$ with
zero first homology, so that $\im\partial_2=\ker\partial_1$, which is exactly
the relation module $\ab{R}$ (cf.~Lemma~\ref{lem:relmodexact}); moreover,
$\ker\partial_2$ is isomorphic to $\pi_2(K)$. Therefore we have an exact
sequence
\begin{equation*}
  \tag{$\dagger$}0\to\pi_2(K)\to(\ZG)^4\xrightarrow{\partial_2}\ab{R}\to0.
\end{equation*}
Applying Schanuel's lemma to the sequences $(*)$ and $(\dagger)$, we deduce
that
\[
\pi_2(K)\oplus(\ZG)^3\cong N\oplus(\ZG)^4.
\]

Now let $L$ be the $2$-complex obtained by taking a $1$-point union of $K$
with $3$ copies of $S^2$. Evidently $\pi_1(L)$ is again isomorphic to
$\Gamma$, and $\pi_2(L)=\pi_2(K)\oplus(\ZG)^3$, which we have seen is
isomorphic to $N\oplus(\ZG)^4$. Attach four $3$-cells to $L$ to fill the four
($\Gamma$-orbits of) $2$-spheres on the right-hand side of this direct sum,
and call the resulting $3$-complex $M$. (The attaching maps of these 3-cells
can be described explicitly by tracing through the above algebra.)
\begin{lemma}
  $M$ enjoys property $D(2)$.
\end{lemma}
\begin{proof}
  Since $\Gamma$ is virtually free, $H^3(\Gamma,\ZG)=0$, and it follows from
  the proof of~\cite[Theorem~3.5]{HarlSome} that $M$ satisfies $D(2)$.
\end{proof}

\begin{proposition}
  Suppose that $\Gamma_{m,n}$ cannot be presented with fewer than
  $4$~relations on the given generators (i.e.~that $\Gamma_{m,n}$ has a
  relation gap). Then $M$ is a counterexample to the $D(2)$~conjecture.
\end{proposition}
\begin{proof}
  By hypothesis, $\defic(\Gamma)=0$ and $\Gamma=\pi_1(M)$. But
  \begin{align*}
    \chi(M)-1&=-4+(4+3)-4\\
    &=-1\\
    &<0\\
    &=\defic(\Gamma),
  \end{align*}
  and so $M$ cannot have the homotopy type of a finite $2$-complex. On the
  other hand, $M$ is a $3$-complex with the $D(2)$~property by the previous
  lemma.
\end{proof}
\noindent For alternative approaches to the $D(2)$~conjecture,
see~\cite{Johnson}.

\section{Concerning lower bounds on deficiency}\label{sec:lower}
\noindent The nub of both the relation gap problem's difficulty and its
attraction is that at present we seem to have no computable invariants that
give lower bounds on $d_F(R)$ without giving an identical bound on the number
of generators that the relation module requires.  There is one
sometimes-computable invariant in the literature that at first sight inspires
hope in regard to the relation gap problem, namely the deficiency test
developed by Martin Lustig in his work on higher Fox ideals~\cite{LustigFox}.
However, in this section we explain why this method cannot help one to
establish relation gaps.  We are grateful to Ian Leary for helpful comments
about this.

Let $\Gamma$ be a group of type $\FPnn{2}$.  Given a free resolution
$(C_*,\partial)$ of $\Z$ over $\ZG$, finitely generated in degrees $\le 2$,
one defines the \emph{second directed Euler characteristic} of $C_*$ to be
$\chid(C_*) = \rk(C_0)-\rk(C_1)+\rk(C_2)$, and one defines $\chid(\Gamma)$ by
taking the infimum over all such resolutions.

Lustig's deficiency test concerns representations of $\ZG$ into non-zero
unital rings $\Lambda$ in which left and right inverses coincide.

Let $K$ be the standard $2$-complex of a finite presentation $\mathscr{P}$ for
$\Gamma$, let $M$ be a $\ZG$-module with $\pi_2(K)\subset M \subset C_2(\tilde
K)$, and identify $C_2(\tilde K)$ with the free $\ZG$-module on the relators
$r_j$ of $\mathscr{P}$ (which index the $\ZG$-orbits of $2$-cells in $\tilde
K$). Fix a finite generating set $\{a_i\}$ for $M$ and express each generator
in terms of the basis of $C_2(\tilde K)$, say $a_i = \sum_j a_{ij}r_j$.
Lustig proves that if there is a unital homomorphism $\rho:\ZG\to\Lambda$ such
that $\rho(a_{ij})=0$ for all $i$ and $j$, then $\mathscr{P}$ realizes the
deficiency of $\Gamma$.

\begin{proposition}\label{prop:lustignogap}
  If there exists a representation $\rho$ satisfying Lustig's criterion, then
  the presentation $\mathscr{P}$ does not have a relation gap.
\end{proposition}
\begin{proof}
  What Lustig actually proves is that $\chid(\Gamma) = \defic(\mathscr{P})$.
  By definition, $\chid(\Gamma)$ is no greater than the directed Euler
  characteristic of the resolution in Lemma~\ref{lem:relmodexact}, which is
  the abelianized deficiency of $\mathscr{P}$. But Lemma~\ref{lem:defdef}
  states that $\adef(\mathscr{P})\le\defic(\mathscr{P})$. Therefore
  $\adef(\mathscr{P})=\defic(\mathscr{P})$.
\end{proof}

\bibliography{books}

\providecommand{\bysame}{\leavevmode\hbox to3em{\hrulefill}\thinspace}
\providecommand{\MR}{\relax\ifhmode\unskip\space\fi MR }
\providecommand{\MRhref}[2]{%
  \href{http://www.ams.org/mathscinet-getitem?mr=#1}{#2}
}
\providecommand{\href}[2]{#2}
\begin{thebibliography}{10}

\bibitem{BestvinaBrady}
M.~Bestvina and N.~Brady, \emph{Morse theory and finiteness properties of
  groups}, Invent. Math. \textbf{129} (1997), 445--470.

\bibitem{BTArtin}
M.~R. Bridson and M.~Tweedale, \emph{Presentations of finite-index subgroups of
  right-angled {A}rtin groups}, in preparation.

\bibitem{EpsteinFinite}
D.~B.~A. Epstein, \emph{Finite presentations of groups and {$3$}-manifolds},
  Quart. J. Math. Oxford Ser. (2) \textbf{12} (1961), 205--212.

\bibitem{GruenbergLinnell}
K.~Gruenberg and P.~Linnell, \emph{Generation gaps and abelianized defects of
  free products}, in preparation.

\bibitem{HarlSome}
J.~Harlander, \emph{Some aspects of efficiency}, Groups---Korea '98 (Pusan), de
  Gruyter, 2000, pp.~165--180.

\bibitem{HausmannWeinberger}
J.-C. Hausmann and S.~Weinberger, \emph{Caract\'eristiques d'{E}uler et groupes
  fondamentaux des vari\'et\'es de dimension {$4$}}, Comment. Math. Helv.
  \textbf{60} (1985), 139--144.

\bibitem{Hillman}
J.~A. Hillman, \emph{An homology {$4$}-sphere group with negative deficiency},
  Enseign. Math. (2) \textbf{48} (2002), 259--262.

\bibitem{HogLustigMetzler}
C.~Hog, M.~Lustig, and W.~Metzler, \emph{Presentation classes, {$3$}-manifolds
  and free products}, Geometry and topology (College Park, Md., 1983/84),
  Lecture Notes in Math., vol. 1167, Springer, 1985, pp.~154--167.

\bibitem{Howie}
J.~Howie, \emph{On locally indicable groups}, Math. Z. \textbf{180} (1982),
  445--461.

\bibitem{HowieBB}
\bysame, \emph{Bestvina--{B}rady groups and the plus construction}, Math.\
  Proc.\ Camb.\ Phil.\ Soc. \textbf{127} (1999), 487--493.

\bibitem{Johnson}
F.~E.~A. Johnson, \emph{Stable modules and the {$D(2)$}-problem}, London
  Mathematical Society Lecture Note Series, vol. 301, Cambridge University
  Press, 2003.

\bibitem{LustigFox}
M.~Lustig, \emph{On the rank, the deficiency and the homological dimension of
  groups: the computation of a lower bound via {F}ox ideals}, Topology and
  combinatorial group theory (Hanover, NH, 1986/1987; Enfield, NH, 1988),
  Lecture Notes in Math., vol. 1440, Springer, 1990, pp.~164--174.

\bibitem{LustigEff}
\bysame, \emph{Non-efficient torsion-free groups exist}, Comm. Algebra
  \textbf{23} (1995), 215--218.

\bibitem{Swan}
R.~G. Swan, \emph{Minimal resolutions for finite groups}, Topology \textbf{4}
  (1965), 193--208.

\bibitem{Wall}
C.~T.~C. Wall, \emph{Finiteness conditions for {${\rm CW}$}-complexes}, Ann. of
  Math. (2) \textbf{81} (1965), 56--69.

\end{thebibliography}
\bibliographystyle{amsplain}

\end{document}